\documentclass[12pt]{article}
\usepackage{amsmath,amsthm,amssymb}
\usepackage{graphicx}

\newcommand{\C}{\mbox{\rm \,l\kern-0.52em C}}
\newcommand{\Ce}{\rm \,l\kern-0.35em C}

\newtheorem{theorem}{Theorem}[section]

\newtheorem{definition}[theorem]{Definition}
\newtheorem{prop}[theorem]{Proposition}

\newtheorem{example}[theorem]{Example}
\newtheorem{lemma}[theorem]{Lemma}
\newtheorem{remark}[theorem]{Remark}

\renewenvironment{proof}{{\bf Proof:}}{\mbox{}\hfill $\Box$}

\theoremstyle{definition}
\date{}

\title{Periodic cyclic homology of crossed products}

\author{Michael Puschnigg}

\begin{document} 



\maketitle

\begin{abstract}
We discuss the cyclic homology of crossed product algebras from the
Cuntz-Quillen point of view. The periodic cyclic homology of a crossed product algebra $A\rtimes G$ 
is described in terms of the $G$-action on periodic cyclic bicomplexes of crossed products of $A$ by the cyclic subgroups of $G$ .
\end{abstract}

\section{Introduction}

The work of Burghelea \cite{Bu} and Nistor \cite{Ni} on the cyclic homology of group rings and crossed product algebras, carried out soon after the invention of cyclic homology fourty years ago, constitutes an integral part of the theory.
We recall their results in a form suitable for our purposes. The key observation is that the commutator quotient 
$$
HH^{A\rtimes G}_0(M,A\rtimes G)=M/[M,A\rtimes G]
\eqno(1.1)
$$ 
of an $A\rtimes G$-bimodule $M$ can be calculated in two steps. At first one forms the commutator quotient $HH_0^{A}(M,A)=M/[M,A]$, which is a $G$-module under the adjoint action, and then one passes to the module of coinvariants under this action. For the corresponding derived functors this implies the existence of a natural isomorphism 
$$
{\mathbb H}_*(G,\,C_*^{A}(M,A))\,\overset{\simeq}{\longrightarrow}\,HH^{A\rtimes G}_*(M,A\rtimes G),
\eqno(1.2)
$$ 
which identifies the Hochschild homology of the pair $(M,A\rtimes G)$ over $A\rtimes G$ with the hyperhomology of $G$ with coefficients in the Hochschild complex of the pair $(M,A)$ over $A$.
For $M=A\rtimes G$ the complex $C_*^{A}(A\rtimes G,A)$ decomposes under the adjoint action into a direct sum labeled by the conjugacy classes of $G$, inducing a similar decomposition for the Hochschild- and cyclic complex of $A\rtimes G$.
For $A=k$ the contribution of the conjugacy classes of the torsion elements to $HC_*(A\rtimes G)$ equals
$$
HC_*(k\rtimes G)_{[\mathrm{tors}]}\,\overset{\simeq}{\longrightarrow}\,\underset{i=0}{\overset{\infty}{\bigoplus}}\,{\mathbb H}_{*-2i}(G,\,\mathrm{Ad}(Vect_k(G_{\mathrm{tors}})),
\eqno(1.3)
$$
the hyperhomology of $G$ with coefficients in the vector space over the set of torsion elements of $G$, equipped with the adjoint action. If $v\in G$ is of infinite order, and if
$$
1\to v^{\mathbb Z}\to \mathrm{Cent}_G(v)\to N(v)\to 1
\eqno(1.4)
$$
is the associated central extension, then Burghelea and Nistor show that the contribution of the conjugacy class $[v]$ to $HC_*(A\rtimes G)$  is a module over                                                                                                                                                 
the cohomology ring $H^*(N(v),k)$ and that Connes' operation 
$$
S:HC_*(A\rtimes G)_{[v]}\to HC_{*-2}(A\rtimes G)_{[v]}
\eqno(1.5)
$$ 
is given by the cap product 
with the extension class $\alpha\in H^2(N(v),k)$ of (1.1).\\
\\
The 40th Birthday Conference of cyclic homology seems to be an appropriate occasion to review the work of Burghelea and Nistor from the point of view of Cuntz and Quillen \cite{CQ1},\cite{CQ2}, which provides a lot of conceptual insight into cyclic theory.
Typical features of the Cuntz-Quillen approach are the strictly ${\mathbb Z}/2{\mathbb Z}$-graded setting, putting the emphasis on periodic cyclic rather than ordinary cyclic or Hochschild theory, and the extensive use of universal algebras.\\
\\
We recall the``derived  functor analogy'' of Cuntz-Quillen, which serves as a guiding principle for our work.
Classical derived functors on the derived category of an abelian category are constructed by replacing an arbitrary chain complex by a well behaved quasi-isomorphic complex and by evaluating the functor to be derived on the latter.
Cuntz and Quillen propose to replace a given algebra by a well behaved topologically nilpotent extension, and to study the given functor, in our case the periodic cyclic bicomplex, on the latter class of algebras.\\
\\
In this spirit we introduce for a given set $X$ a natural nilpotent extension $k\langle X\rangle$ of the ground field $k$. In particular, if a group $G$ acts on $X$, then $k\langle X\rangle$ becomes a $G$-algebra. For $X=G$, equipped with the translation action, and a $G$-algebra $A$ the $G$-algebra $A\langle G\rangle\,=\,A\otimes k\langle G\rangle$ is then
a natural nilpotent extension of $A$ which is free as a $G$-module. It serves as our ``well behaved model'' of $A$ and will play a key role throughout this paper.\\
\\
In order to formulate our results conveniently we introduce periodic cyclic bicomplexes $\widehat{CC}(A:R)$, given by the usual direct product of tensor powers of the given algebra $A$, but this time relative to a not necessarily commutative ground ring $R$.\\
\\
The cyclic bicomplex $\widehat{CC}_*(A\langle G\rangle\rtimes G:k\rtimes G)_{[e]}$ has then two canonical bases and can be written in two ways. On the one hand we find the bicomplex
$\widehat{CC}_*(A\rtimes G)_{[e]}$, and on the other hand the space of $G$-coinvariants of the complex $\widehat{CC}_*(A\langle G\rangle)$ of free $G$-modules. So $HP_*(A\rtimes G)_{[e]}$ equals the hyperhomology 
of $G$ with coefficients in $\widehat{CC}_*(A\langle G\rangle)$. Goodwillie's theorem on the invariance of periodic cyclic homology under nilpotent extensions allows to identify this with the hyperhomology of $G$ with coefficients in $\widehat{CC}_*(A)$, recovering the result of Nistor \cite{Ni}. A similar picture emerges for the contributions of the other conjugacy classes. 
In terms of the central extension (1.1), we find for elements $v\in G$ of finite order a natural isomorphism
$$
HP_*(A\rtimes G)_{[v]}\,\overset{\simeq}{\longrightarrow}\,\widehat{{\mathbb H}}_*\left(N_v,\,\widehat{CC}(A\rtimes v^{\mathbb Z}:k\rtimes v^{\mathbb Z})_{\{v\}}\right).
\eqno(1.6)
$$

It follows that the elliptic part of the periodic cyclic homology of a crossed product, i.e. the contribution of the conjugacy classes of all torsion elements except the unit,  is given by 
$$
HP_*(A\rtimes G)_{\mathrm{ell}}\,\simeq\,\widehat{{\mathbb H}}_*\left(G,\,\underset{\underset{1<\vert v\vert<\infty}{v\in G}}{\widehat{\bigoplus}}\,\widehat{CC}_*(A\rtimes v^{{\mathbb Z}}: k\rtimes v^{{\mathbb Z}})_{\{v\}}\right).
\eqno(1.7)
$$
For conjugacy classes of elements $v\in G$ of infinite order we obtain partial results similar to those of Nistor \cite{Ni}. 
If $v$ acts trivially on $A$ a formula similar to (1.6) continues to hold with the right hand side replaced by its derived inverse limit with respect to the iterated cap product with the cohomology class of the extension (1.4).\\
\\
Roughly speaking one could say that the periodic cyclic homology of a crossed product $A\rtimes G$ is determined by a family of cyclic bicomplexes of the crossed products of $A$ by the cyclic subgroups of $G$ and the $G$-action on this family.\\
\\
Our approach does not only provide quite explicit and rather simple formulas for the contributions of various conjugacy classes to the periodic cyclic (co)homology of a crossed product, but extends also nicely to a topological setting. The attempt of adapting Nistor's original approach to Banach crossed products was, on the contrary, a quite frustrating experience for the author and motivated the present research.\\
\\
I thank the referee for a very thorough reading of the manuscript, his suggestions for improvement, and for pointing out an error in the first version of the paper.

\section{Preliminaries}
\subsection{Base change of cyclic complexes}
Fix a field $k$ of characteristic 0 and a unital associative $k$-algebra $R$.
\begin{definition}
An $R$-algebra is a unitary $R$-bimodule $A$, together with an $R$-bimodule homomorphism  
$$
m:\,A\otimes_R A\,\to\,A,
\eqno(2.1)
$$
which defines a (not necessarily unital) $k$-algebra structure on $A$.
A homomorphism $\varphi:A\to B$ of $R$-algebras is a map which is simultaneously a homomorphism of $R$-bimodules and of $k$-algebras.
\end{definition}
With this definition $R$-algebras form a category. 
 
\begin{definition} 
Let $A$ be an $R$-algebra. The graded $R$-bimodule $$\Omega^*(A:R)\,=\,\underset{n=0}{\overset{\infty}{\bigoplus}}\,\Omega^n(A:R)
\eqno(2.2)
$$ of
 {\bf algebraic differential forms} of $A$ over $R$ is defined as
$$
\begin{array}{ccc}
\Omega^n(A:R) & \simeq & A^{\otimes_R^{n+1}} \oplus A^{\otimes_R^{n}} \\
& &  \\
a_0da_1\ldots da_n & \leftrightarrow & a_0\otimes a_1\otimes\ldots\otimes a_n  \\
& &  \\
da_1\ldots da_n & \leftrightarrow &  a_1\otimes\ldots\otimes a_n  \\
\end{array}
\eqno(2.3)
$$
with $R$-bimodule-structure induced by the natural one on the tensor powers of $A$. 
The commutator quotient with respect to this bimodule structure is the graded $k$-vector space 
$$
\Omega^*(A:R)_\natural\,=\,\Omega^*(A:R)/[\Omega^*(A:R),R].
\eqno(2.4)
$$
\end{definition}
Note that this differs from the notion of relative differential forms in \cite{CQ1}. 
Under the canonical linear projection 
$$
\Omega^*A=\Omega^*(A:k)_\natural\,\to\,\Omega^*(A:R)_\natural
\eqno(2.5)
$$
 the well known Hochschild- and Connes-operators 
$$
\begin{array}{cc}
b:\,\Omega^*A\,\to\,\Omega^{*-1}A, & B:\,\Omega^*A\,\to\,\Omega^{*+1}A
\end{array}
\eqno(2.6)
$$
descend to linear operators 
$$
\begin{array}{cc}
b:\,\Omega^*(A:R)_\natural\,\to\,\Omega^{*-1}(A:R)_\natural, & B:\,\Omega^*(A:R)_\natural\,\to\,\Omega^{*+1}(A:R)_\natural,
\end{array}
\eqno(2.7)
$$
satisfying the usual identities
$$
b^2\,=\,B^2\,=\,bB+Bb\,=\,0.
\eqno(2.8)
$$
\begin{definition}
\item[i)] The {\bf cyclic bicomplex } of $A$ over $R$ is given by 
$$
CC_*(A:R)\,=\,\left(\Omega^*(A:R)_\natural,\,\,b+B\right).
\eqno(2.9)
$$
This is a ${\mathbb Z}/2{\mathbb Z}$-graded naturally contractible chain complex. 
\item[ii)] The {\bf Hodge filtration} of the cyclic bicomplex is the descending filtration by the subcomplexes 
$$
\mathrm{Fil}^m_{\mathrm{Hodge}}CC_*(A:R)=\underset{n=m}{\overset{\infty}{\bigoplus}}\Omega^n(A:R)_\natural\oplus b(\Omega^m(A:R)_\natural),\,m\in{\mathbb N},
\eqno(2.10)
$$
generated by the algebraic differential forms of degree at least $m$. 
\item[iii)] The associated graded complex is quasi-isomorphic to the Hochschild complex
$$
C_*(A:R)\,=\,\left(\Omega^*(A:R)_\natural,\,\,b\right).
\eqno(2.11)
$$
\item[iv)] The {\bf periodic cyclic bicomplex} of $A$ over $R$ is the completion of $CC_*(A:R)$ with respect to the topology defined by the Hodge filtration:
$$
\widehat{CC}_*(A:R)\,=\,\left( \underset{n=0}{\overset{\infty}{\prod}}\,\Omega^n(A:R)_\natural,\,\,b+B\right).
\eqno(2.12)
$$
This is a complete ${\mathbb Z}/2{\mathbb Z}$-graded chain complex, the grading being given by the parity of forms.\\
\item[v)] The {\bf periodic cyclic homology} $HP_*(A:R)$ of $A$ over $R$ is the homology of the periodic cyclic bicomplex $\widehat{CC}_*(A:R)$.
The {\bf periodic cyclic cohomology} $HP^*(A:R)$ of $A$ over $R$ is the cohomology of the topologically dual complex of $k$-linear functionals on $\widehat{CC}_*(A:R)$ which are continuous (i.e., vanishing on  
$\mathrm{Fil}^m_{\mathrm{Hodge}}\widehat{CC}_*(A:R)$ for $m\gg 0$).
\end{definition}

Every homomorphism $f:A\to B$ of $R$-algebras gives rise to a continuous morphism $\widehat{CC}(f):\,\widehat{CC}_*(A:R)\,\to\,\widehat{CC}_*(B:R)$ of periodic cyclic bicomplexes. 
So the periodic cyclic bicomplex defines a functor from the category of $R$-algebras to the category of ${\mathbb Z}/2{\mathbb Z}$-graded 
topologically complete chain complexes. \\
If $S\subset R$ is a unital subalgebra, the canonical projection $\Omega^*(A:S)_\natural\,\to\,\Omega^*(A:R)_\natural$ induces a morphism
$$
\pi_{S,R}:\,\widehat{CC}_*(A:S)\,\longrightarrow\,\widehat{CC}_*(A:R)
\eqno(2.13)
$$
of topologically complete chain complexes.

\begin{example}
For unital $A$ one has
$$
\begin{array}{ccc}
\widehat{CC}(A:A)\,=\,\underset{n=0}{\overset{\infty}{\prod}}\,A/[A,A], & HP_0(A:A)=A/[A,A], & HP_1(A:A)=0.\\
\end{array}
\eqno(2.14)
$$
\end{example}

\subsection{Auxiliary algebras}

Let $X$ be a set and let $k\langle X\rangle$ be the $k$-vector space over $X$. The projection
$$
p:\,X\times X\to X,\,\,(x,y)\mapsto y
\eqno(2.15)
$$
turns $k\langle X\rangle$ into an associative $k$-algebra with multiplication
$$
m:\,k\langle X\rangle\otimes_kk\langle X\rangle\,=\,k
\langle X\times X\rangle\,\overset{k\langle p\rangle}{\longrightarrow}\,k\langle X\rangle.
\eqno(2.16)
$$
If $X$ consists of a single element the algebra thus obtained is canonically isomorphic to the ground field $k$. Every map of sets $f:X\to Y$ gives rise to an algebra homomorphism 
$$
k\langle f\rangle:\,k\langle X\rangle\to k\langle Y\rangle.
\eqno(2.17)
$$
In particular, every action of a group $G$ on $X$ induces a $G$-action on $k\langle X\rangle$ by $k$-algebra automorphisms.
The constant map to a point defines an augmentation homomorphism
$$
\epsilon_X:\,k\langle X\rangle\to k.
\eqno(2.18)
$$
The multiplication in $k\langle X\rangle$ is characterized by the identity
$$
a\cdot b\,=\,\epsilon_X(a)b,\,\,\,\forall a,b\in k\langle X\rangle.
\eqno(2.19)
$$
In particular, the augmentation ideal $I_\epsilon\,=\,\mathrm{Ker}(\epsilon_X)$ satisfies 
$$
I_\epsilon^2\,=\,I_\epsilon\cdot k\langle X\rangle\,=\,0.
\eqno(2.20)
$$
The tensor product of a unital $k$-algebra $A$ and $k\langle X\rangle$ in the category of $k$-algebras is denoted by
$$
A\langle X\rangle = A\otimes_kk\langle X\rangle.
\eqno(2.21)
$$
The algebra extension 
$$
0\,\longrightarrow\,I_\epsilon\otimes_k A\,\longrightarrow\,A\langle X\rangle\,\overset{\epsilon_A}{\longrightarrow}\,A\,\longrightarrow\,0
\eqno(2.22)
$$
is nilpotent by (2.20) and splits (non-canonically) as extension of $k$-algebras. If a group $G$ acts on the $k$-algebra $A$ and acts also freely on the set $X$, the diagonal action turns 
$A\langle X\rangle$ into a $G$-algebra whose underlying $G$-module is free.

\section{Cyclic complexes of crossed products}

\subsection{Cyclic complexes attached to crossed products}
Let $G$ be a group and let $A$ be a unital $G$-algebra over $k$. We let $k\langle G\rangle$ carry the $G$-action induced by left translation and 
equip $A\langle G\rangle$ with the diagonal $G$-action. 

\begin{lemma}
 Let $U\subset G$ be a subgroup and let $\sigma:U\backslash G\to G$ be a set theoretic section of the canonical projection $G\to U\backslash G$.
\begin{itemize}
\item[i)] The linear map 
$$
\begin{array}{cccc}
\iota_\sigma: & (A\rtimes U)\langle U\backslash G\rangle & \to & A\langle G\rangle\rtimes U, \\
 & & & \\
 & (u_g a)\langle x\rangle & \mapsto & u_g (a\langle \sigma(x)\rangle) \\
 \end{array}
 \eqno(3.1)
$$
is a homomorphism of $k$-algebras.
\item[ii)] The composition of chain maps
$$
\begin{array}{ccc}
CC_*((A\rtimes U)\langle U\backslash G\rangle) & \overset{CC(\iota_\sigma)}{\longrightarrow} & CC_*(A\langle G\rangle\rtimes U) \\
&& \\
\parallel & & \downarrow\pi_{k,k\rtimes U} \\
& & \\
CC_*((A\rtimes U)\langle U\backslash G\rangle) & \overset{\simeq}{\longrightarrow} & CC_*(A\langle G\rangle\rtimes U:\,k\rtimes U)\\
\end{array}
\eqno(3.2)
$$
is an isomorphism of filtered bicomplexes (with respect to Hodge filtrations).
\end{itemize}
\end{lemma}

\begin{proof}
i) The map $\iota_\sigma$ is an algebra homomorphism as
$$
\iota_\sigma((u_ga)\langle x\rangle)\iota_\sigma((u_hb)\langle y\rangle) =u_ga\langle \sigma(x)\rangle u_hb\langle \sigma(y)\rangle=u_{gh}h^{-1}(a\langle \sigma(x)\rangle)b\langle \sigma(y)\rangle=
$$
$$
=u_{gh}h^{-1}(a)\langle h^{-1}\sigma(x)\rangle b\langle \sigma(y)\rangle=u_{gh}h^{-1}(a)b\langle h^{-1}\sigma(x)\rangle\langle \sigma(y)\rangle=u_{gh}h^{-1}(a)b\langle\sigma(y)\rangle=
$$
$$=\iota_\sigma((u_{gh}h^{-1}(a)b)\langle y\rangle)=
\iota_\sigma((u_ga)\langle x\rangle\cdot (u_hb)\langle y\rangle).
$$
ii) Note that the $k$-linear map
$$
\begin{array}{ccc}
kU\otimes_k A\langle U\backslash G\rangle\otimes _k kU & \longrightarrow & A\langle G\rangle\rtimes U \\
 & & \\
 u_g\otimes a\langle x\rangle\otimes u_h & \mapsto & u_ga\langle \sigma(x)\rangle u_h=u_{gh}h^{-1}(a)\langle h^{-1}\sigma(x)\rangle \\
\end{array}
\eqno(3.3)
$$
is an isomorphism of $kU$-bimodules, so that $A\langle G\rangle\rtimes U$ becomes a free $kU$-bimodule with basis $A\langle U\backslash G\rangle$. Consequently
$$
 (A\langle G\rangle\rtimes U)^{\otimes_{kU}^n}\simeq(kU\otimes_k A\langle U\backslash G\rangle\otimes _k kU)^{\otimes_{kU}^n}\simeq(kU\otimes_k A\langle U\backslash G\rangle)^{\otimes_{k}^n}\otimes_k kU
$$
and
$$
 (A\langle G\rangle\rtimes U)^{\otimes_{kU}^n}/[ (A\langle G\rangle\rtimes U)^{\otimes_{kU}^n},kU]\,\simeq\,(kU\otimes_k A\langle U\backslash G\rangle)^{\otimes_{k}^n}\,\simeq\,((A\rtimes U)\langle U\backslash G\rangle)^{\otimes_{k}^n}
$$
as $k$-vector spaces, which shows that the composition of the two chain maps is an isomorphism. The second assertion follows from the fact that all occuring
maps preserve the degree of differential forms.
\end{proof}

In the sequel we denote by $\mathrm{Ad}(U)$ the set $U$ with the adjoint $U$-action, and let $\mathrm{Vect}_k(\mathrm{Ad}(U))$ be the associated  $k$-vector space with linear $U$-action.

\begin{lemma}
Let $i:A\langle G\rangle\to A\langle G\rangle\rtimes U$ be the canonical inclusion. The linear map
$$
\begin{array}{ccc}
\mathrm{Vect}_k(\mathrm{Ad}(U))\otimes_k\Omega^*(A\langle G\rangle) & \longrightarrow & \Omega^*(A\langle G\rangle\rtimes U) \\
 & & \\
 u_g\otimes\omega & \mapsto & u_g\cdot i_*(\omega)
 \end{array}
\eqno(3.4)
$$
fits into a commutative diagram
$$
\begin{array}{ccc}
\mathrm{Vect}_k(\mathrm{Ad}(U))\otimes_k\Omega^*(A\langle G\rangle) & \longrightarrow & \Omega^*(A\langle G\rangle\rtimes U) \\
 & & \\
 \downarrow && \downarrow \\
 & & \\
 \left(\mathrm{Vect}_k(\mathrm{Ad}(U))\otimes_k\Omega^*(A\langle G\rangle)\right)_U & \overset{\simeq}{\longrightarrow} & \Omega^*(A\langle G\rangle\rtimes U:k\rtimes U)_\natural \\
 \end{array}
 \eqno(3.5)
$$
whose lower horizontal arrow is an isomorphism. 
Here the left vertical arrow is the projection onto the space of coinvariants under the diagonal $U$-action and the right vertical arrow is given by the canonical projection of absolute onto relative algebraic differential forms.
Every subgroup of $G$ centralizing $U$ acts naturally on this diagram.
\end{lemma}

\begin{proof}
Let $\Phi$ be the composition of the upper horizontal and the right vertical map. We find for $h\in U$ 
$$
\Phi(h(u_{g}\otimes a_0\langle g_0\rangle\otimes\ldots\otimes a_n\langle g_n\rangle))=\Phi(u_{hgh^{-1}}\otimes h(a_0\langle g_0\rangle)\otimes\ldots\otimes h(a_n\langle g_n\rangle))
$$

$$
=u_{h}u_{g}u_{h^{-1}}h(a_0\langle g_0\rangle)\otimes\ldots\otimes h(a_n\langle g_n\rangle)=u_{g}u_{h^{-1}}h(a_0\langle g_0\rangle)\otimes\ldots\otimes h(a_n\langle g_n\rangle)u_h
$$

$$
=u_{g}u_{h^{-1}}h(a_0\langle g_0\rangle) u_h\otimes\ldots\otimes u_{h^{-1}}h(a_n\langle g_n\rangle)u_h=\Phi(u_{g}\otimes a_0\langle g_0\rangle\otimes\ldots\otimes a_n\langle g_n\rangle),
$$
which shows that $\Phi$ is constant on $U$-orbits. 
Using the identification of\\ $\Omega^*((A\langle G\rangle)\rtimes U:k\rtimes U) _\natural$ with $\Omega^*((A\rtimes G)\langle U\backslash G\rangle)$ of Lemma 3.1, one may construct a $k$-linear section $\sigma'$ of $\Phi$ by putting
$$
\sigma'(u_{g_0}a_0\langle x_0\rangle\otimes\ldots\otimes u_{g_{n-1}}a_{n-1}\langle x_{n-1}\rangle\otimes u_{g_n}a_n\langle x_n\rangle)\,=
\eqno(3.6)
$$

$$
u_{g_0\ldots g_n}\otimes (g_1\ldots g_n)^{-1}(a_0\langle \sigma(x_0)\rangle)\otimes\ldots\otimes g_n^{-1}(a_{n-1}\langle\sigma(x_{n-1})\rangle)\otimes a_n\langle\sigma(x_n)\rangle.
$$

The image of $\sigma'$ intersects each $U$-orbit in exactly one element. This proves our assertion.
\end{proof}

\subsection{The homogeneous decomposition}

We recall the well known homogeneous decomposition of cyclic complexes of crossed product algebras \cite{Bu},\cite{Ni}.
Let $G$ be a group and denote by $[v]$ the conjugacy class of $v\in G$. For a $G$-algebra $A$ over $k$ let
$\Omega^*(A\rtimes G)_{[v]}$ be the linear span of the differential forms
$$
a_0u_{g_0}d(a_1u_{g_1})\ldots d(a_nu_{g_n})\in\Omega^*(A\rtimes G), g_0g_1\ldots g_n\in[v],
$$
and
$$
d(b_1u_{h_1})\ldots d(b_nu_{h_n})\in\Omega^*(A\rtimes G),h_1\ldots h_n\in[v].
$$
Put
$$
CC_*(A\rtimes G)_{[v]}\,=\,(\Omega^*(A\rtimes G)_{[v]},b+B).
\eqno(3.7)
$$
There are canonical decompositions
$$
\begin{array}{ccc}
\Omega^*(A\rtimes G)\,=\,\underset{[v]}{\bigoplus}\,\Omega^*(A\rtimes G)_{[v]} & \text{and} & CC_*(A\rtimes G)\,=\,\underset{[v]}{\bigoplus}\,CC_*(A\rtimes G)_{[v]},
\end{array}
\eqno(3.8)
$$
where the sum runs over the set of conjugacy classes of $G$. If $U\subset G$ is a subgroup, the homogeneous decompositions (3.8) descend under the canonical projections (2.13) to similar decompositions 
$$
\Omega^*(A\rtimes G:k\rtimes U)_\natural\,=\,\underset{[v]}{\bigoplus}\,\Omega^*(A\rtimes G:k\rtimes U)_{\natural\, [v]} 
\eqno(3.9)
$$
and
$$ 
CC_*(A\rtimes G:k\rtimes U)\,=\,\underset{[v]}{\bigoplus}\,CC_*(A\rtimes G:k\rtimes U)_{[v]}.
\eqno(3.10)
$$
The corresponding completions $\widehat{CC}_*(A\rtimes G)_{[v]}$ and $\widehat{CC}_*(A\rtimes G:k\rtimes U)_{[v]}$ are topological direct summands of $\widehat{CC}_*(A\rtimes G)$ and $\widehat{CC}_*(A\rtimes G:k\rtimes U)$, respectively, but the latter complexes are in general
neither the direct sum nor the direct product of their factors in the homogeneous decomposition.
The partition of the set of conjugacy classes of $G$ into the singleton given by the conjugacy class of the unit, the union of the other conjugacy classes of elements of finite order, and the union of the conjugacy classes of elements of infinite order 
gives rise to a decomposition of the cyclic bicomplex, its completion, and the periodic cyclic (co)homology of a crossed product algebra into the direct sum of three factors, called their {\bf homogeneous}, {\bf elliptic} and {\bf inhomogeneous} parts.\\
\\For $v\in G$ let $Z_v$ be its centralizer in $G$. The cyclic subgroup $U$ of $G$, generated by $v$, fits then into a central extension
$$
1\,\longrightarrow\,U=v^{{\mathbb Z}}\,\longrightarrow\,Z_v\,\longrightarrow\,N_v\,\longrightarrow\,1.
\eqno(3.11)
$$

\begin{prop}
The canonical chain map 
$$
CC_*(A\langle G\rangle\rtimes U:\,k\rtimes  U)\,\longrightarrow\,CC_*(A\langle G\rangle\rtimes G:\,k\rtimes  G)\,\overset{(3.2)}{\simeq}\,CC_*(A\rtimes G)
$$
induces an isomorphism
$$
\left(CC_*(A\langle G\rangle\rtimes U:\,k\rtimes  U)_{\{v\}}\right)_{Z_v}\,\,\,\overset{\simeq}{\longrightarrow}\,\,\,CC_*(A\rtimes G)_{[v]}
\eqno(3.12)
$$
of filtered chain complexes (with respect to the Hodge filtrations).
\end{prop}

\begin{proof}
It is clear that the chain map preserves Hodge filtrations. Therefore it suffices to show that the underlying map of vector spaces is an isomorphism.
By Lemma 3.2, applied to $U=v^{{\mathbb Z}}$, and (3.9) there is a canonical linear isomorphism
$$
\Omega^*(A\langle G\rangle)_U\,\overset{\simeq}{\longrightarrow}\,\Omega^*(A\langle G\rangle\rtimes U:\,k\rtimes  U)_{\natural\,\{v\}}
\eqno(3.13)
$$
By construction it commutes with the action of the centralizer $Z_v$, so that one obtains an isomorphism of the spaces of $Z_v$-coinvariants,
 which fits into the commutative diagram
 $$
\begin{array}{ccc}
\Omega^*(A\langle G\rangle)_{Z_v} & \overset{\simeq}{\longrightarrow} & \left(\Omega^*(A\langle G\rangle\rtimes U:\,k\rtimes U)_{\natural\,\{v\}}\right)_{Z_v} \\
 & & \\
 \downarrow & & \downarrow \\
 & & \\
\left(Vect_k([v])\otimes_k\Omega^*(A\langle G\rangle)\right)_{G} & \overset{\simeq}{\longrightarrow} & \Omega^*(A\langle G\rangle\rtimes G:\,k\rtimes G)_{\natural\,[v]}\\
 \end{array}
 \eqno(3.14)
 $$
 The upper horizontal arrow is an isomorphism by (3.13), the lower one by Lemma 3.2, applied to $U=G$, and the left vertical arrow is an isomorphism because $Z_v$ is the stabilizer of the element $v\in[v]$ under the transitive adjoint action of $G$ on $[v]$.
 It results that the right vertical arrow is an isomorphism, which is our assertion.
\end{proof}

\subsection{Hyperhomology}

It is well known and easy to prove that quasi-isomorphic ${\mathbb Z}_+$-graded chain complexes of $G$-modules have isomorphic hyperhomology groups ${\mathbb H}_*(G,-)$. In general this is no longer true for ${\mathbb Z}/2{\mathbb Z}$-graded complexes.
This forces us to add this technical section dealing with several cases where the previous result continues to hold.\\
\\
We consider the category of ${\mathbb Z}/2{\mathbb Z}$-graded complexes of $k\rtimes G$-modules, equipped with an adic topology, defined by a decreasing chain of subcomplexes. 
The morphisms are given by continuous, $G$-equivariant chain maps. A morphism which becomes an isomorphism in the associated chain-homotopy category is called a continuous equivariant chain homotopy equivalence.
A filtration defining the given topology will be called admissible.
The completion of such a chain complex in the sense of general topology equals
$$
\widehat{C}_*\,\simeq\,\underset{\underset{n}{\longleftarrow}}{\lim}\,C_*/\mathrm{Fil}^nC_*
\eqno(3.15)
$$
for any admissible filtration. 
The algebraic tensor product $C_*\otimes_k C'_*$ of two complexes is topologized by
$$
\mathrm{Fil}_n(C_*\otimes_k C'_*)=\underset{p+q=n}{\sum}\,\mathrm{Fil}_pC_*\otimes_k \mathrm{Fil}_qC'_*
\eqno(3.16)
$$
This is independent of the choice of the admissible filtrations on the individual factors. 
The completion $\widehat{C_*}\widehat{\otimes}_k\widehat{C'_*}$ of the tensor product $C_*\otimes_k C'_*$  contains in general the algebraic tensor product $\widehat{C_*}\otimes_k\widehat{C'_*}$ of the individual completions as dense subcomplex. The complex of $G$-coinvariants of $C_*$ equals
$$
(C_*)_G\,=\,C_*\otimes_{k\rtimes G}k
\eqno(3.17)
$$
equipped with the filtration induced by the given one on $C_*$ and the constant filtration on $k$.
It is the largest quotient complex of $C_*$ on which $G$ acts trivially. 
The complex $(C_*)_G$ is complete if $C_*$ is. 
The Hom-complex $L(C_*,C'_*)$ is the complex of all $k$-linear maps from $C_*$ to $C'_*$, equipped with the differential
$$
\partial_{L(C_*,C'_*)}(\varphi)\,=\,\partial_{C'_*}\circ\varphi-(-1)^{\mathrm{deg}(\varphi)}\varphi\circ\partial_{C_*},\,\,\,\varphi\in L(C_*,C'_*).
\eqno(3.18)
$$
 Its subcomplex  of continuous linear maps is denoted by ${\mathcal L}(C_*,C'_*)$.
\\
To a ${\mathbb Z}_+$-graded chain complex $P_*$ one may associate the ${\mathbb Z}/2{\mathbb Z}$-graded chain complex
$$
(P^{\pm}_*,\partial),\,P^{\pm}_{\mathrm{ev}}=\underset{n=0}{\overset{\infty}{\bigoplus}}\,P_{2n},\,P^{\pm}_{\mathrm{odd}}=\underset{n=0}{\overset{\infty}{\bigoplus}}\,P_{2n+1},
\eqno(3.19)
$$
topologized by the filtration
$$
\mathrm{Fil}^mP^{\pm}_*\,=\,\partial P_m\oplus\underset{n=m}{\overset{\infty}{\bigoplus}}P_n
\eqno(3.20)
$$
Its completion equals
$$
\widehat{P}^\pm_*\,=\,\underset{n=0}{\overset{\infty}{\prod}}\,P_n.
\eqno(3.21)
$$

\begin{definition}
Let $G$ be a group and let $P_*$ be a projective resolution of the constant $G$-module $k$. Let $C_*$ be a  ${\mathbb Z}/2{\mathbb Z}$-graded chain complex
of $k\rtimes G$-modules, topologized by a decreasing family of subcomplexes, and let $C^*={\mathcal L}(C_*,k)={\mathcal L}(\widehat{C}_*,k)={\widehat C}^*$ be the topologically dual complex.
\begin{itemize}
\item[i)] The {\bf hyperhomology} of $G$ with coefficients in $C_*$  equals
$$
{\mathbb H}_*(G,\,C_*)\,=\,H_*((P^\pm_*\otimes_k C_*)_G).
\eqno(3.22)
$$
\item[ii)] The {\bf continuous hyperhomology} of $G$ with coefficients in $\widehat{C_*}$ equals
$$
\widehat{{\mathbb H}}_*(G,\,\widehat{C}_*)\,=\,H_*((\widehat{P}^\pm_*\widehat{\otimes}_k \widehat{C}_*)_G).
\eqno(3.23)
$$
\item[iii)] The {\bf hypercohomology} of $G$ with coefficients in $C_*$  equals
$$
{\mathbb H}^*(G,\,C_*)\,=\,H^*(L(P^\pm_*\otimes_k C_*,k)^G).
\eqno(3.24)
$$
\item[iv)] The {\bf continuous hypercohomology} of $G$ with coefficients in $\widehat{C}^*$ equals
$$
\widehat{{\mathbb H}}^*(G,\,\widehat{C}_*)\,=\,H^*({\mathcal L}(\widehat{P}^\pm_*\widehat{\otimes}_k\widehat{C}_*, k)^G).
\eqno(3.25)
$$

\end{itemize}
\end{definition}
This is (up to canonical isomorphism) independent of the choice of resolution. The continuous hyper(co)homology depends only on the topology of $\widehat{C}_*$, but not on the particular choice of filtration used to define it.
 If there exists a resolution of finite length of the constant $G$-module $k$ by finitely generated, projective $kG$-modules, then the natural maps
 $$
 \begin{array}{ccc}
 {\mathbb H}_*(G,\,{\widehat C_*})\to \widehat{{\mathbb H}}_*(G,\,{\widehat C_*}) & \text{and} & \widehat{{\mathbb H}}^*(G,\,{\widehat C_*})\to {\mathbb H}^*(G,\,{\widehat C_*})\\
 \end{array}
 \eqno(3.26)
 $$ 
are isomorphisms.

\begin{lemma}
Let $A$ be an $R$-algebra and let $G$ be a group acting on $A$ by $R$-algebra automorphisms such that $\Omega^*(A:R)_\natural$ becomes a degree-wise free $kG$-module. Let  $P_*$ be a projective resolution of the constant $G$-module $k$. Then the natural chain map
$$
\epsilon\otimes_k\mathrm{id}:\,(P^\pm_*\otimes CC_*(A:R))_G\,\longrightarrow\,(CC_*(A:R))_G
\eqno(3.27)
$$
becomes a continuous chain-homotopy equivalence after completion. In particular
$$
\widehat{{\mathbb H}}_*(G,\widehat{CC}_*(A:R))\,\overset{\simeq}{\longrightarrow}\,H_*(\widehat{CC}_*(A:R)_G).
\eqno(3.28)
$$
\end{lemma}

\begin{proof}
Let $h:\widetilde{P}_*\to\widetilde{P}_{*+1},*\geq -1,$ be a $k$-linear contracting chain homotopy of the augmented resolution $\widetilde{P}_*$ of the constant $G$-module $k$. By assumption
$\Omega^n(A:R)_\natural$ is a free $G$-module i.e., $\Omega^n(A:R)_\natural\simeq kG\otimes_k V_n$ as $kG$-module for some $k$-vector space $V_n$ and all $n\geq 0$.
Let $\widetilde{h}:\widetilde{P}^\pm_*\otimes_k CC_*(A:R)\to \widetilde{P}^\pm_{*+1}\otimes_k CC_*(A:R)$ be the unique morphism of $kG$-modules satisfying $\widetilde{h}\vert_{\widetilde{P}_*\otimes_k V_*}\,=\,h\otimes_k\mathrm{id}$.
The chain map 
$$\varphi=\mathrm{id}-(\partial\circ\widetilde{h}+\widetilde{h}\circ\partial):\widetilde{P}^\pm_*\otimes_k CC_*(A:R)\to \widetilde{P}^\pm_*\otimes_k CC_*(A:R)$$ satisfies 
$$
\varphi(\widetilde{P}_{\geq n}\otimes_k\Omega^{\geq m}(A:R)_\natural)\,\subset\,\widetilde{P}_{\geq n+1}\otimes_k\Omega^{\geq m-1}(A:R)_\natural.
\eqno(3.29)
$$
So the infinite series 
$$
\nu\,=\,\underset{j=0}{\overset{\infty}{\sum}}\,\widetilde{h}\circ\varphi^j.
\eqno(3.30)
$$
converges to a continuous $G$-equivariant contracting chain homotopy of the completion of  $\widetilde{P}^\pm_*\otimes_k CC_*(A:R)$. It follows that the augmentation morphism
$$
\epsilon\otimes_k\mathrm{id}:\,(P^\pm_*\otimes CC_*(A:R))_G\,\longrightarrow\,(CC_*(A:R))_G
\eqno(3.31)
$$
becomes a chain homotopy equivalence after completion. Taking homology one obtains a natural isomorphism
$$
\widehat{{\mathbb H}}_*(G,\,\widehat{CC}_*(A:R))\,\overset{\simeq}{\longrightarrow}\,H_*(\widehat{CC}_*(A:R)_{G}).
$$
\end{proof}

A similar proof establishes the following lemma for groups of finite cohomological dimension over $k$. As we do not want to make this restriction a more elaborate argument is necessary.

\begin{lemma}
Let $A$ be a unital $G$-algebra over $k$ and let $U\subset G$ be a normal subgroup. Let $X$ be a $G$-set such that $\epsilon_A:\,A\langle X\rangle\,\to\,A$ has a $U$-equivariant $k$-algebra section.
Then 
$$
\epsilon_{A_*}:\,\widehat{{\mathbb H}}_*(G,\widehat{CC}_*(A\langle X\rangle\rtimes U:kU))\,\overset{\simeq}{\longrightarrow}\,\widehat{{\mathbb H}}_*(G,\widehat{CC}_*(A\rtimes U:kU)).
\eqno(3.32)
$$
is a natural isomorphism compatible with homogeneous decompositions.
\end{lemma}

\begin{proof}
Let $\varphi:A\to A\langle X\rangle$ be a $U$-equivariant homomophism of $k$-algebras that splits the canonical epimorphism $\epsilon_A:\,A\langle X\rangle\to A$.  These homomorphisms
 induce chain maps 
$$
(\epsilon_A\rtimes\mathrm{id})_*:\,\widehat{CC}(A\langle X\rangle\rtimes U:k\rtimes U)\,\longrightarrow\,\widehat{CC}(A\rtimes U:k\rtimes U)
$$
and
$$
(\varphi\rtimes\mathrm{id})_*:\,\widehat{CC}(A\rtimes U:k\rtimes U)\,\longrightarrow\,\widehat{CC}(A\langle X\rangle\rtimes U:k\rtimes U).
$$
Moreover
$$
s\mapsto (1-s)(\varphi\circ\epsilon_A)\rtimes\mathrm{id}+s\cdot\mathrm{id},\,s\in[0,1],
\eqno(3.33)
$$
defines an affine homotopy of $k\rtimes U$-algebra homomorphisms between $(\varphi\circ\epsilon_A)\rtimes\mathrm{id}: A\langle X\rangle\rtimes U\to A\langle X\rangle\rtimes U$ and the identity. The corresponding chain homotopy of periodic cyclic complexes  (see [Lo], page 117) between $\widehat{CC}((\varphi\circ\epsilon_A)\rtimes\mathrm{id})$ and the identity descends to a continuous linear Cartan homotopy operator
$$
h_0:\,\widehat{CC}_*(A\langle X\rangle\rtimes U:k\rtimes U)\,\to\,\widehat{CC}_{*+1}(A\langle X\rangle\rtimes U:k\rtimes U).
\eqno(3.34)
$$
Note that
$$
h_0((A\langle X\rangle\rtimes U)\,\Omega^n(A\langle X\rangle\rtimes U))\subset\Omega^{n+1}(A\langle X\rangle\rtimes U)
\eqno(3.35)
$$
by (2.19).
We cannot use this operator directly if $U$ is of infinite cohomological dimension because it does not preserve Hodge filtrations. To overcome this difficulty we use an auxiliary preliminary homotopy.
Recall that the cone of a morphism $f:C_*\to C'_*$ of complexes is defined as the complex 
$$
\mathrm{Cone}(f)_*\,=\,C_{*}\oplus C'_{*+1},\,\,\,\partial_{\mathrm{Cone}}=
\left(
\begin{matrix}
\partial_{C_*} & 0 \\
f & -\partial_{C'_*} \\
\end{matrix}
\right).
\eqno(3.36)
$$
Put 
$$
h_1\,=\,\left(
\begin{matrix}
h_0 & \widehat{CC}(\varphi\rtimes\mathrm{id}) \\
0 & 0 \\
\end{matrix}
\right):\,\mathrm{Cone}_*(\widehat{CC}(\epsilon_A\rtimes\mathrm{id}))\to \mathrm{Cone}_{*+1}(\widehat{CC}(\epsilon_A\rtimes\mathrm{id})).
\eqno(3.37)
$$
This is a $k$-linear contracting homotopy of $\mathrm{Cone}_*(\widehat{CC}(\epsilon_A\rtimes\mathrm{id}))$. \\
Let $h_2:\mathrm{Cone}_*(\widehat{CC}(\epsilon_A\rtimes\mathrm{id}))\to \mathrm{Cone}_{*+1}(\widehat{CC}(\epsilon_A\rtimes\mathrm{id}))$ be the $k$-linear map sending 
$\omega\in d\Omega^*(A\langle X\rangle\rtimes U)$ to $-d(\varphi(1_A) u_e)\omega$
and annihilating $(A\langle X\rangle\rtimes U)d\Omega^*(A\langle X\rangle\rtimes U)$, and sending $\omega'\in d\Omega^*(A\rtimes U)$ to $d(1_Au_e)\omega'$ and annihilating $(A\rtimes U)d\Omega^*(A\rtimes U).$
Let 
$$
\varphi_2:=\mathrm{id}-(h_2\circ\partial_{\mathrm{Cone}(\widehat{CC}(\epsilon_A\rtimes \mathrm{id}))}+\partial_{\mathrm{Cone}(\widehat{CC}(\epsilon_A\rtimes\mathrm{id}))}\circ h_2)
$$ 
and put 
$$
h_3=h_1\circ\varphi_2+h_2:\mathrm{Cone}_*(\widehat{CC}(\epsilon_A\rtimes\mathrm{id}))\to \mathrm{Cone}_{*+1}(\widehat{CC}(\epsilon_A\rtimes\mathrm{id})).
\eqno(3.38)
$$
The operator $h_3$ is again a $k$-linear contracting homotopy of $\mathrm{Cone}_*(\widehat{CC}(\epsilon_A\rtimes\mathrm{id}))$. But in addition, it strictly increases the degree of algebraic differential forms.

Let $P_*$ be a degree-wise free resolution i.e., $P_n\simeq kG\otimes_k W_n,n\geq 0,$ of the constant $G$-module $k$ (it always exists).  We let
$$
h:\,P^\pm_*\otimes_k \mathrm{Cone}_*(\widehat{CC}(\epsilon_A\rtimes\mathrm{id}))\to P^\pm_*\otimes_k \mathrm{Cone}_{*+1}(\widehat{CC}(\epsilon_A\rtimes\mathrm{id}))
\eqno(3.39)
$$
be the unique $G$-equivariant linear operator such that $h_*\vert_{W_*\otimes_k \mathrm{Cone}}=\mathrm{id}\otimes_k h_3$ and put 
$$
\psi=\mathrm{id}-(\partial\circ h+h\circ\partial):\,P^\pm_*\otimes_k \mathrm{Cone}(\widehat{CC}(\epsilon_A\rtimes\mathrm{id}))_*\to P^\pm_*\otimes_k \mathrm{Cone}(\widehat{CC}(\epsilon_A\rtimes\mathrm{id}))_*.
$$
The operator $\underset{n=0}{\overset{\infty}{\sum}}\,h\circ\psi^n$
 is then a well defined, continuous and $G$-equivariant contracting chain homotopy of $\widehat{P}^\pm_*\widehat{\otimes}_k \mathrm{Cone}_*(\widehat{CC}(\epsilon_A\rtimes\mathrm{id})).$
Consequently 
$$
\mathrm{id}\otimes_k \widehat{CC}(\epsilon_A\rtimes\mathrm{id}):\,\widehat{P}^\pm_*\widehat{\otimes}_k\widehat{CC}_*(A\langle G\rangle\rtimes U:k\rtimes U)\,\to\,\widehat{P}^\pm_*\widehat{\otimes}_k\widehat{CC}_*(A\rtimes U: k\rtimes U)
\eqno(3.40)
$$
is a continuous $G$-equivariant chain homotopy equivalence and
$$
(\epsilon_{A}\rtimes\mathrm{id})_*:\,\widehat{{\mathbb H}}_*(G,\widehat{CC}_*(A\langle G\rangle\rtimes U:k\rtimes U))\,\overset{\simeq}{\longrightarrow}\,\widehat{{\mathbb H}}_*(G,\widehat{CC}_*(A\rtimes U:k\rtimes U))
$$
is an isomorphism. 
\end{proof}

\section{Periodic cyclic homology of crossed products}

\begin{theorem}
Let $G$ be a group and let $A$ be a $G$-algebra over $k$. Let $v\in G$ and let 
$$
1\,\longrightarrow\,U=v^{{\mathbb Z}}\,\longrightarrow\,Z_v\,\longrightarrow\,N_v\,\longrightarrow\,1
$$
be the associated central extension. Then there is a natural isomorphism
$$
HP_*(A\rtimes G)_{[v]}\,\overset{\simeq}{\longrightarrow}\,\widehat{{\mathbb H}}_*\left(N_v,\widehat{CC}_*(A\langle Z_v\rangle\rtimes U:k\rtimes U)_{\{v\}}\right).
\eqno(4.1)
$$
\end{theorem}

\begin{proof}
By Proposition 3.3 there is an isomorphism
$$
\left(CC_*(A\langle G\rangle\rtimes U:\,k\rtimes  U)_{\{v\}}\right)_{Z_v}\,\,\,\overset{\simeq}{\longrightarrow}\,\,\,CC_*(A\rtimes G)_{[v]}
$$
of filtered chain complexes. By definition  (see Lemma 3.2) the restriction of the canonical $Z_v$-action on $\Omega^*(A\langle G\rangle\rtimes U:\,k\rtimes U)_{\natural}$ to $U$ is trivial
so that 
$$
HP_*(A\rtimes G)_{[v]}\,\overset{\simeq}{\longrightarrow}\,H_*((\widehat{CC}_*(A\langle G\rangle\rtimes U:k\rtimes U)_{\{v\}})_{N_v}).
$$
From (3.13) we learn that it is a degree-wise free $k\rtimes N_v$-module. So by Lemma 3.5 there is a natural isomorphism
$$
\widehat{{\mathbb H}}_*\left(N_v,\widehat{CC}_*(A\langle G\rangle\rtimes U:k\rtimes U)_{\{v\}}\right)\,\overset{\simeq}{\longrightarrow}\,H_*((\widehat{CC}_*(A\langle G\rangle\rtimes U:k\rtimes U)_{\{v\}})_{N_v}).
$$
Lemma 3.6 may be applied to both projections in the diagram
$$
A\langle G\rangle \leftarrow A\langle G\times Z_v\rangle \to A\langle Z_v\rangle
$$
and obtains the natural quasi-isomorphism
$$
\widehat{{\mathbb H}}_*\left(N_v,\widehat{CC}_*(A\langle Z_v\rangle\rtimes U:k\rtimes U)_{\{v\}}\right) \,\overset{\simeq}{\longrightarrow}\,  
\widehat{{\mathbb H}}_*\left(N_v,\widehat{CC}_*(A\langle G\rangle\rtimes U:k\rtimes U)_{\{v\}}\right)
$$
This concludes the proof of the theorem.
\end{proof}

\subsection{The homogeneous part}

\begin{theorem}
Let $G$ be a group and let $A$ be a $k$-algebra on which $G$ acts by automorphisms. There are natural isomorphisms
$$
HP_*(A\rtimes G)_{[e]}\,\overset{\simeq}{\longrightarrow}\,\widehat{{\mathbb H}}_*(G,\,\widehat{CC}_*(A)),
\eqno(4.2)
$$
and
$$
HP^*(A\rtimes G)_{[e]}\,\overset{\simeq}{\longrightarrow}\,\widehat{{\mathbb H}}^*(G,\,\widehat{CC}_*(A)),
\eqno(4.3)
$$
which identify the homogeneous part of the periodic cyclic (co)homology of the crossed product $A\rtimes G$ with
the continuous hyper(co)homology of $G$ with coefficients in the periodic cyclic bicomplex of $A$.
\end{theorem}

\begin{proof}
We apply Theorem 4.1 to the case $v=e,\,U=1,\,N_v=G$ and obtain a natural isomorphism
$$
HP_*(A\rtimes G)_{[e]}\,\overset{\simeq}{\longrightarrow}\,\widehat{{\mathbb H}}_*(G,\,\widehat{CC}_*(A\langle G\rangle)).
$$
Goodwillie's Theorem in the strengthened form of Lemma 3.6 shows that 
$$
\widehat{{\mathbb H}}_*(G,\,\widehat{CC}_*(A\langle G\rangle))\,\overset{\simeq}{\longrightarrow}\,\widehat{{\mathbb H}}_*(G,\,\widehat{CC}_*(A)),
$$
is a natural isomorphism as well. A similar argument establishes (4.3).
\end{proof}

The {\bf bivariant periodic cyclic cohomology} \cite{CQ2} of a pair $(A,B)$ of\\ $k$-algebras is defined as the homology of the ${\mathbb Z}/2{\mathbb Z}$-graded chain complex 
$$
\widehat{CC}_*(A,B)\,=\,{\mathcal L}(\widehat{CC}_*(A),\widehat{CC}_*(B)).
\eqno(4.4)
$$
We recall that the assumptions of the following two theorems are satisfied if some classifying space of $G$ has the homotopy type of a finite $CW$-complex.

\begin{theorem}
Let $G$ be a group such that the constant $G$-module $k$ possesses a resolution of finite length by finitely generated free $kG$-modules. Let $A$ be a $G$-algebra over $k$ and let $B$ be a $k$-algebra.
Then there is a natural isomorphism
$$
HP_*(A\rtimes G,B)_{[e]}\,\overset{\simeq}{\longrightarrow}\,{\mathbb H}^*(G,\,\widehat{CC}_*(A,B))
\eqno(4.5)
$$
identifying the homogeneous part of the bivariant periodic cyclic cohomology of the pair $(A\rtimes G,B)$ with the hypercohomology of $G$ with coefficients in the bivariant periodic cyclic bicomplex of the pair $(A,B)$.
\end{theorem}

\begin{proof}

We learn from Proposition 3.3, applied to  $v=e$, that there is a canonical isomorphism
$$
(CC_*(A\langle G\rangle))_G\,\overset{\simeq}{\longrightarrow}\,\,\,CC_*(A\rtimes G)_{[e]}
$$
of chain complexes preserving Hodge filtrations.
Let $P_*$ be a resolution of finite length of the constant $G$-module $k$ by finitely generated free $kG$-modules.
As $\Omega^*(A\langle G\rangle)$
is degree-wise free as $kG$-module 
$$
 (P^\pm_*\otimes_k CC_*(A\langle G\rangle))_G\,\longrightarrow\,(CC_*(A\langle G\rangle))_G
$$
becomes a continuous chain-homotopy equivalence after completion by Lemma 3.5. 
Lemma 3.6, and in particular assertion (3.40), applied to $U=\{e\},\,X=G,$ and\\ $\varphi:A\to A\langle G\rangle, a\mapsto a\langle e\rangle$, show that the left hand arrow in the natural diagram
$$
(P^\pm_*\otimes_k CC_*(A))_G\,\longleftarrow\, (P^\pm_*\otimes_k CC_*(A\langle G\rangle))_G\,\longrightarrow\,CC_*(A\rtimes G)_{[e]}
\eqno(4.6)
$$
becomes a continuous chain homotopy equivalence after completion, too. \\
There are natural isomorphisms of chain complexes
$$
\mathrm{Hom}_{kG}\left(\widehat{P}^\pm_*,\,{\mathcal L}(\widehat{CC}_*(A),\widehat{CC}_*(B))\right)\,\overset{\simeq}{\longrightarrow}\,
\mathrm{Hom}_k\left(\widehat{P}^\pm_*,\,{\mathcal L}(\widehat{CC}_*(A),\widehat{CC}_*(B))\right)^G
$$

$$
\overset{\simeq}{\longrightarrow}\,
{\mathcal L}\left(\widehat{P}^\pm_*\widehat{\otimes}_k \widehat{CC}_*(A),\widehat{CC}_*(B)\right)^G
\overset{\simeq}{\longrightarrow}\,
{\mathcal L}\left((\widehat{P}^\pm_*\widehat{\otimes}_k \widehat{CC}_*(A))_G,\widehat{CC}_*(B)\right).
$$

Making use of (4.6) again and taking cohomology we arrive at the desired result.
\end{proof}

\begin{theorem}
Let $G$ be a group such that the constant $G$-module $k$ possesses a resolution of finite length by finitely generated free $kG$-modules. Let $A$ be a $k$-algebra and let $B$ be a $G$-algebra over $k$.
Then there is a natural isomorphism
$$
HP_*(A,B\rtimes G)_{[e]}\,\overset{\simeq}{\longrightarrow}\,{\mathbb H}_*(G,\,\widehat{CC}_*(A,B))
\eqno(4.7)
$$
identifying the homogeneous part of the bivariant periodic cyclic cohomology of the pair $(A,B\rtimes G)$ with the hyperhomology of $G$ with coefficients in the bivariant periodic cyclic bicomplex of the pair $(A,B)$.
\end{theorem}

\begin{proof}
Let $P_*$ be a resolution of finite length of the constant $G$-module $k$ by finitely generated free $kG$-modules. So $P_*\simeq kG\otimes_k W_*$ as graded $kG$-modules for some finite dimensional graded $k$-vector space $W_*$.
For a finite dimensional vector space $V$ we denote by $\check{V}$ its dual. Let $\widehat{({P}^\pm_*\otimes _k CC_*(B))_G}$ be the completion of $({P}^\pm_*\otimes _k CC_*(B))_G$ with respect to the filtration topology. 
There are natural isomorphisms of linear spaces 
$$
{\mathcal L}\left(\widehat{CC}_*(A),\widehat{({P}^\pm_*\otimes _k CC_*(B))_G}\right)\,\overset{\simeq}{\longrightarrow}\,
{\mathcal L}\left(\widehat{CC}_*(A),W^\pm_*\otimes _k\widehat{CC}_*(B)\right)
$$

$$
\overset{\simeq}{\longrightarrow}\,
{\mathcal L}\left(\widehat{CC}_*(A),{\check{\check W}}^\pm_*\otimes _k\widehat{CC}_*(B)\right)\,
\overset{\simeq}{\longrightarrow}\,{\mathcal L}\left(\widehat{CC}_*(A),{\mathcal L}(\check{W}^\pm_*,\widehat{CC}_*(B))\right)
$$

$$
\overset{\simeq}{\longrightarrow}\,{\mathcal L}\left(\check{W}^\pm_*\otimes_k\widehat{CC}_*(A),\widehat{CC}_*(B)\right)\,
\overset{\simeq}{\longrightarrow}\,Hom_k\left(\check{W}^\pm_*,{\mathcal L}(\widehat{CC}_*(A),\widehat{CC}_*(B))\right)
$$

$$
\overset{\simeq}{\longrightarrow}\,\check{\check{W}}_*^{\pm}\otimes_k{\mathcal L}(\widehat{CC}_*(A),\widehat{CC}_*(B))\,
\overset{\simeq}{\longrightarrow}\,\left(P_*^{\pm}\otimes_k{\mathcal L}(\widehat{CC}_*(A),\widehat{CC}_*(B))\right)_G.
$$

It is easily verified that this is a morphism of chain complexes. Making use of (4.6) again and taking homology, we prove our claim as before.
\end{proof}

\subsection{The elliptic part}

\begin{theorem}
Let  $G$ be a group and let $A$ be a $k$-algebra on which $G$ acts by automorphisms. 
There is a natural isomorphism
$$
HP_*(A\rtimes G)_{\mathrm{ell}}\,\simeq\,\widehat{{\mathbb H}}_*\left(G,\,\underset{\underset{1<\vert v\vert<\infty}{v\in G}}{\widehat{\bigoplus}}\,\widehat{CC}_*(A\rtimes U: k\rtimes U)_{\{v\}}\right)
\eqno(4.8)
$$
where the direct sum is labeled by the set of torsion elements in $G$ and we denote by $U$ the finite cyclic subgroup generated by the torsion element $v$. 
The subscripts $\{v\}$ and $[v]$ denote the conjugacy classes of $v$ in $U$ and $G$, respectively. The complex $\underset{\underset{1<\vert v\vert<\infty}{v\in G}}{\widehat{\bigoplus}}\,\widehat{CC}(A\rtimes U: k\rtimes U)_{\{v\}}$ is the completion of 
$\underset{\underset{1<\vert v\vert<\infty}{v\in G}}{\bigoplus}\,CC_*(A\rtimes U: k\rtimes U)_{\{v\}},$ equipped with the direct sum of the individual filtrations.
For a single torsion element one obtains a natural isomorphism
$$
HP_*(A\rtimes G)_{[v]}\,\simeq\,\widehat{{\mathbb H}}_*(Z_v,\,\widehat{CC}_*(A\rtimes U:\,k\rtimes U)_{\{v\}})
\eqno(4.9)
$$
where $Z_v$ is the centralizer of $v$ in $G$.
\end{theorem}

The theorem expresses the elliptic part of the periodic cyclic homology of the crossed product $A\rtimes G$ of $A$ by $G$ in terms 
of the hyperhomology of $G$ with values in the completed direct sum of the periodic cyclic bicomplexes of the crossed products of $A$ by the various finite cyclic subgroups of $G$.

\begin{proof}
Let $v\in G$ be an element of finite order and let $U\subset G$ be the finite cyclic group generated by $v$. By Proposition 3.3 there is an isomorphism of filtered complexes
$$
\left(CC_*(A\langle G\rangle\rtimes U:k\rtimes U)_{\{v\}}\right)_{Z_v}\,\,\,\overset{\simeq}{\longrightarrow}\,\,\,CC_*(A\rtimes G)_{[v]}.
$$
As $k$ is of characteristic zero every projective $k\rtimes N_v$-module is projective as $k\rtimes Z_v$-module as well.
Let $P_*$ be a projective resolution of the $G$-module $k$. It is at the same time a projective resolution of the $Z_v$-module $k$.  
The finite central subgroup $U\subset Z_v$ acts trivially on the complex $CC_*(A\langle G\rangle\rtimes U:\,k\rtimes U)_{\{v\}}$ The underlying vector space 
$\Omega^*(A\langle G\rangle\rtimes U:k\rtimes U)_\natural$ is a degree-wise free $N_v$-module and thus projective as $Z_v$-module.
Lemma 3.5 implies therefore that 
$$
\left(P^\pm_*\otimes_k CC_*(A\langle G\rangle\rtimes U:k\rtimes U)_{\{v\}}\right)_{Z_v}\,
\longrightarrow\,\left(CC_*(A\langle G\rangle\rtimes U:k\rtimes U)_{\{v\}}\right)_{Z_v}
$$
becomes a continuous chain homotopy equivalence after completion.  Consider now the canonical homomorphism 
$\epsilon_A:\,A\langle G\rangle\,\longrightarrow\,A.$ 
Its linear section
$$
\varphi:\,A\longrightarrow A\langle G\rangle,\,a\,\mapsto\,\frac{1}{\vert U\vert}\underset{h\in U}{\sum}\,a\langle h\rangle
$$
is actually a homomorphism of $k\rtimes U$-algebras. So one may deduce from Lemma 3.6 and (3.40) that
$$
\left(P^\pm_*\otimes_k CC_*(A\langle G\rangle\rtimes U:k\rtimes U)_{\{v\}}\right)_{Z_v}\,
\longrightarrow\,\left(P^\pm_*\otimes_k CC_*(A\rtimes U:k\rtimes U)_{\{v\}}\right)_{Z_v}
$$
becomes a continuous chain homotopy equivalence after completion. Thus one derives from the previously constructed chain maps after completion and passage to homology the natural diagram of isomorphisms
$$
\widehat{{\mathbb H}}_*(Z_v,\widehat{CC}_*(A\rtimes U:k\rtimes U)_{\{v\}})\,\overset{\simeq}{\leftarrow}\,\widehat{{\mathbb H}}_*(Z_v,\widehat{CC}_*(A\langle G\rangle\rtimes U:k\rtimes U)_{\{v\}})\,\overset{\simeq}{\to}\,HP_*(A\rtimes G)_{[v]},
$$
which is our second claim. Our previous arguments show that there exists a natural chain map
$$
\left(\underset{w\in[v]}{\bigoplus}\,P^\pm_*\otimes_k CC_*(A\langle G\rangle\rtimes w^{{\mathbb Z}}:k\rtimes w^{{\mathbb Z}})_{\{w\}}\right)_G\,\longrightarrow\,CC_*(A\rtimes G)_{[v]}
$$
which becomes a continuous chain homotopy equivalence after completion. Passing to direct sums over all conjugacy classes of torsion elements we obtain a chain map 
$$
\left(\underset{\underset{1<\vert v\vert<\infty}{v\in G}}{\bigoplus}\,P^\pm_*\otimes_k CC_*(A\langle G\rangle\rtimes U:k\rtimes U)_{\{v\}}\right)_G\,\longrightarrow\,CC_*(A\rtimes G)_{\mathrm{ell}}
$$
with similar properties because the filtration shifts are the same in all direct summands. 
Passing to completions and taking homology, we arrive at the first assertion of the theorem.
\end{proof}

There are corresponding results for the elliptic part of periodic cyclic cohomology and bivariant periodic cyclic homology of crossed products. The details are left to the reader.

\subsection{The inhomogeneous part}

In general it seems to be difficult to express the inhomogeneous part of the periodic cyclic homology of a crossed product in terms of simpler homological invariants.
We consider here only a special case treated already by Nistor \cite{Ni}.
\\
Recall that the extension
$$
1\,\longrightarrow\,U=v^{{\mathbb Z}}={\mathbb Z}\,\longrightarrow\,Z_v\,\overset{\pi}{\longrightarrow}\,N_v\,\longrightarrow\,1.
\eqno(4.10)
$$

is classified up to isomorphism by its extension class 
$$
\alpha\in H^2(N_v,{\mathbb Z}).
\eqno(4.11)
$$
If $\sigma:N_v\to Z_v$ is a set theoretic section of $\pi$, 
a homogeneous group cocycle representing $\alpha$ is given by
$$
c_\alpha:\,(N_v)^3\,\to\,{\mathbb Z},\,\,\,(h_0,h_1,h_2)\,\mapsto\,\sigma(h_1^{-1}h_2)\sigma(h_0^{-1}h_2)^{-1}\sigma(h_0^{-1}h_1).
\eqno(4.12)
$$

We follow the classical approch and begin with the calcuation of the inhomogeneous part of Hochschild homology. 
For an $A$-bimodule $M$ let
$$
H_*(C_*(M,A))=HH_*(M,A)=Tor^{A\otimes_k A^{op}}_*(M,A)
\eqno(4.13)
$$ 
be the Hochschild homology of $(M,A)$. Let $G$ be a group acting on $A$.
For $v\in G$ let $A^{(v)}$ be the $A$-bimodule with underlying vector space $A$ and bimodule structure
$$
A\otimes_k A^{(v)}\otimes_k A\to A^{(v)},\,a'\otimes a\otimes a''\mapsto v^{-1}(a')aa''.
\eqno(4.14)
$$

\begin{theorem}
Let $G$ be a group and let $A$ be a unital $G$-algebra. Let $v\in G$ be an element of infinite order. There is a natural isomorphism
$$
HH_*(A\rtimes G)_{[v]}\,\overset{\simeq}{\longrightarrow}\,{\mathbb H}_*\left(N_v,C_*(A^{(v)},A)\otimes^{{\mathbb L}}_{k\rtimes U}k\right).
\eqno(4.15)
$$
\end{theorem}

\begin{proof}
We proceed in several steps.\\
Step 1:\\
Lemma 3.1, Lemma 3.2 and Proposition 3.3 hold for Hochschild complexes as well as for cyclic complexes. So there is an natural isomorphism
$$
HH_*(A\rtimes G)_{[v]}\,\overset{\simeq}{\longrightarrow}\,H_*(C_*(A\langle G\rangle^{(v)},A\langle G\rangle)_{Z_v}).
\eqno(4.16)
$$
Step 2:\\
The Eilenberg-Zilber theorem in Hochschild theory \cite{Lo} states that for unital $k$-algebras $A,B$, $A$-bimodules $M$ and $B$-bimodules $N$ there is a natural chain-homotopy equivalence
$$
C_*(M\otimes_kN,A\otimes_k B)\,\overset{\Delta}{\longrightarrow}\,C_*(M,A)\otimes_k C_*(N,B).
\eqno(4.17)
$$
Step 3:\\
The algebra $k\langle X\rangle$ is not unital, but naturally $H$-unital \cite{Lo} in the sense that
$$
h':\,C'_*(k\langle X\rangle)\,\to\,C'_*(k\langle X\rangle),\,\langle x_0\rangle\otimes\ldots\otimes\langle x_n\rangle\mapsto\langle x_0\rangle\otimes\ldots\otimes\langle x_n\rangle\otimes\langle x_n\rangle
\eqno(4.18)
$$
is a contracting chain-homotopy of the $b'$-complex of $k\langle X\rangle$ which is natural with respect to maps $X\to Y$ of sets. It follows then from step 2 and excision in Hochschild homology \cite{Lo} that
 there is a $Z_v$-equivariant chain homotopy equivalence
$$
C_*((A\langle G\rangle)^{(v)},A\langle G\rangle)\,\overset{\Delta}{\longrightarrow}\,C_*(A^{(v)},A)\otimes_k C_*(k\langle G\rangle).
$$
(Note that $k\langle G\rangle^{(v)}=k\langle G\rangle$ by (2.19).) So one may pass to coinvariants and obtains a natural chain-homotopy equivalence
$$
\left(C_*((A\langle G\rangle)^{(v)},A\langle G\rangle)\right)_{Z_v}\,\overset{\Delta}{\longrightarrow}\,\left(C_*(A^{(v)},A)\otimes_k C_*(k\langle G\rangle)\right)_{Z_v}.
\eqno(4.19)
$$
Step 4:\\
Observe that the natural $G$-equivariant chain map
$$
\begin{array}{ccc}
C_*^{\mathrm{Bar}}(G,k) & \hookrightarrow & (\Omega^*(k\langle G\rangle),b) \\
 & & \\
 {[}g_0,\ldots,g_n{]} & \mapsto & \langle g_0\rangle d\langle g_1\rangle\ldots d\langle g_n\rangle \\
 \end{array}
 \eqno(4.20)
$$
is equivariantly linearly split with equivariantly contractible cokernel. In particular there is a natural chain-homotopy equivalence
$$
\left(C_*(A^{(v)},A)\otimes_k C_*(k\langle G\rangle)\right)_{Z_v}\,\overset{\sim}{\longrightarrow}\,\left(C_*(A^{(v)},A)\otimes_k C^{Bar}_*(G,k)\right)_{Z_v}.
\eqno(4.21)
$$
As the Bar-complex of a group is a free resolution of the constant $G$-module one deduces from (4.16), (4.19), and (4.21) a natural isomorphism
$$
HH_*(A\rtimes G)_{[v]}\,\overset{\simeq}{\longrightarrow}\,{\mathbb H}_*\left(Z_v,C_*(A^{(v)},A)\right).
$$
Step 5:\\
Observe finally the isomorphism
$$
-\otimes^{{\mathbb L}}_{k\rtimes Z_v}k\,\simeq\, (-\otimes^{{\mathbb L}}_{k\rtimes U}k)\otimes^{{\mathbb L}}_{k\rtimes N_v}k
$$
of derived functors from ${\mathbb D}_-(k\rtimes Z_v)$ to ${\mathbb D}_-(k)$, which implies
$$
{\mathbb H}_*\left(Z_v,C_*(A^{(v)},A)\right)\,\simeq\,H_*\left(C_*(A^{(v)},A)\otimes^{{\mathbb L}}_{k\rtimes Z_v}k\right)\,\simeq
$$

$$
H_*\left((C_*(A^{(v)},A)\otimes^{{\mathbb L}}_{k\rtimes U}k)\otimes^{{\mathbb L}}_{k\rtimes N_v}k\right)\,\simeq\,{\mathbb H}_*\left(N_v,C_*(A^{(v)},A)\otimes^{{\mathbb L}}_{k\rtimes U}k\right).
$$
The theorem is proved.
\end{proof}

\begin{prop}
Let $A$ be a $G$-algebra and let $v\in G$ be an element of infinite order. Then 
$HP_*(A\rtimes G)_{[v]}$ is a $H^*(N_v,k)$-module such that the extension class $\alpha\in H^2(N_v,k)$ of $\mathrm{(4.10)}$ acts as the identity.
\end{prop}

\begin{proof}
We proceed in several steps.\\
Step 1:\\
A weak version of the Eilenberg-Zilber theorem in periodic cyclic homology \cite{Pu}, valid for non-unital algebras as well, states the existence of a natural chain map
$$
\Delta_{cyc}:\,\widehat{CC}(A\otimes_k B)\,\longrightarrow\,\widehat{CC}(A)\widehat{\otimes}_k\widehat{CC}(B)
\eqno(4.22)
$$
for every pair $A,B$ of $k$-algebras, where on the right hand side the completion of the algebraic tensor product of the periodic cyclic complexes of $A$ and $B$ occurs. 
In particular, if $A$ and $B$ are $G$-algebras, then $\Delta$ is $G$-equivariant equivariant. It is compatible with base change in the sense that it descends under the natural projections 
to a natural chain map 
$$
\Delta_{cyc}:\,\widehat{CC}(A\otimes_kB:R\otimes_k S)\,\longrightarrow\,\widehat{CC}(A:R)\,\widehat{\otimes}_k\,\widehat{CC}(B:S).
\eqno(4.23)
$$
This is a natural chain homotopy equivalence if $A$ and $B$ are unital.\\
Step 2:\\
In step 4 of the proof of the previous theorem we observed the close connection between the Bar-complex $C_*^{Bar}(G,k)$ of a group $G$ with coefficients in $k$ and the Hochschild complex of $k\langle G\rangle$. 
In fact, every alternating homogeneous group cocycle $c_\beta\in Z^n(G,k)$ defines a $G$-invariant periodic cyclic cocycle $c'_\beta$ on $k\langle G\rangle$ by the formula
$$
c'_\beta(\langle g_0\rangle d\langle g_1\rangle\ldots d\langle g_n\rangle)=c_\beta([g_0,\ldots,g_n]),\,\,\,c'_\beta(d\langle g_1\rangle\ldots d\langle g_n\rangle)=0,\,\,\, g_0,\ldots,g_n\in G.
\eqno(4.24)
$$
We apply this in the case $G=N_v$.
By abuse of language we call the $N_v$-equivariant natural chain map 
$$
\cap_{c_\beta}:\,\widehat{CC}_*(A\langle Z_v\rangle\rtimes U:k\rtimes U)_{\{v\}}\,\overset{(\mathrm{id}\times\pi)\circ\mathrm{diag}}{\longrightarrow}\,\widehat{CC}_*(A\langle Z_v\times N_v\rangle\rtimes U:k\rtimes (U\times 1))_{\{(v,1)\}}
$$

$$
\overset{\Delta_{cyc}}{\longrightarrow}\,\widehat{CC}_*(A\langle Z_v\rangle\rtimes U:k\rtimes U)_{\{v\}}
\widehat{\otimes}_k\widehat{CC}_*(k\langle N_v\rangle)\,\overset{\mathrm{id}\otimes c'_\beta}{\longrightarrow}\,\widehat{CC}_*(A\langle Z_v\rangle\rtimes U:k\rtimes U)_{\{v\}}
\eqno(4.25)
$$
the ``cap-product'' with the cocycle $c_\beta$. Up to chain homotopy it only depends on the cohomology class of $c_\beta$.
Via the identification (4.1) it turns $HP_*(A\rtimes G)_{[v]}$ into a module over the cohomology ring $H^*(N_v,k)$ of $N_v$ with coefficients in $k$.\\
Step 3: \\
Let $j:U\to(k,+)$ be the homomorphism sending the generator $v\in U$ to 1. Let  $\eta:\,\widehat{CC}(k\langle Z_v\rangle\rtimes U:kU)_{\{v\}}\to k$ be the $N_v$-invariant continuous $k$-linear functional given on one-forms by
$$
\eta(u_v\langle g_0\rangle d\langle g_1\rangle)\,=\,\frac12\left(j(g_0\sigma(g_0^{-1}g_1)g_1^{-1})-j(g_1\sigma(g_1^{-1}g_0)g_0^{-1}) \right),\,\,\,\eta(u_vd\langle g_1\rangle)\,=\,1
\eqno(4.26)
$$
and vanishing on forms of other degrees. The $N_v$-equivariant continuous linear operator
$$
\widetilde{\eta}:\,\widehat{CC}_*(A\langle Z_v\rangle\rtimes U:k\rtimes U)_{\{v\}}\,\overset{\mathrm{diag}}{\longrightarrow}\,\widehat{CC}_*(A\langle Z_v\times Z_v\rangle\rtimes (U\times U):k\rtimes(U\times U))_{\{(v,v)\}}
$$
$$
\overset{\Delta_{cyc}}{\longrightarrow}\,\widehat{CC}_*(A\langle Z_v\rangle\rtimes U:kU)_{\{v\}}\widehat{\otimes}_k \widehat{CC}_*(k\langle Z_v\rangle\rtimes U:kU)_{\{v\}}
$$
$$
\overset{\mathrm{id}\otimes\eta}{\longrightarrow}\,
\widehat{CC}_{*+1}(A\langle Z_v\rangle\rtimes U:k\rtimes U)_{\{v\}}
\eqno(4.27)
$$
defines then an $N_v$-equivariant chain homotopy between $\cap c_{\alpha}$ and the identity, where $c_\alpha$ is the cocycle (4.12) representing the extension class (4.11).
\end{proof}

\begin{theorem}
Let  $G$ be a group and let $A$ be a $k$-algebra on which $G$ acts by automorphisms. Let $v\in G$ be an element of infinite order that acts trivially on $A$ and let $[v]$ be its conjugacy class.
There is a natural isomorphism
$$
HP_*(A\rtimes G)_{[v]}\,\simeq\,{\widehat{\mathbb H}}_*(N_v,\widehat{CC}_*(A))[\alpha^{-1}]
\eqno(4.28)
$$
where
$$
{\widehat{\mathbb H}}_*(N_v,\widehat{CC}_*(A))[\alpha^{-1}]\,=\,{\mathbb R}\underset{\underset{\cap\alpha}{\longleftarrow}}{\lim}\,{\widehat{\mathbb H}}_*(N_v,\widehat{CC}_*(A))
\eqno(4.29)
$$
fits into the natural exact sequence 
$$
0\,\to\,\underset{\underset{\cap\alpha}{\leftarrow}}{\lim}^1({\widehat{\mathbb H}}_{*-1}(N_v,\widehat{CC}_*(A))\,\to
\,{\widehat{\mathbb H}}_*(N_v,\widehat{CC}_*(A))[\alpha^{-1}]
$$
$$
\to\,\underset{\underset{\cap\alpha}{\leftarrow}}{\lim}({\widehat{\mathbb H}}_{*}(N_v,\widehat{CC}_*(A))
\,\to\,0.
\eqno(4.30)
$$
Here the derived projective limit is taken over the iterated cap product 
$$
\cap\alpha:\,{\widehat{\mathbb H}}_*(N_v,\widehat{CC}_*(A))\,\longrightarrow\,{\widehat{\mathbb H}}_*(N_v,\widehat{CC}_*(A))
\eqno(4.31)
$$ 
with the extension class of $\mathrm{(4.10)}$.
\end{theorem}

\newpage
\begin{theorem}
Let  $G$ be a group and let $A$ be a $k$-algebra on which $G$ acts by automorphisms. Let $v\in G$ be an element of infinite order that acts trivially on $A$ and let $[v]$ be its conjugacy class.
There is a natural isomorphism
$$
HP^*(A\rtimes G)_{[v]}\,\simeq\,{\widehat{\mathbb H}}^*(N_v,\widehat{CC}_*(A))[\alpha^{-1}]
\eqno(4.32)
$$
where
$$
{\widehat{\mathbb H}}^*(N_v,\widehat{CC}_*(A))[\alpha^{-1}]\,=\,\underset{\underset{\cup\alpha}{\longrightarrow}}{\lim}\,{\widehat{\mathbb H}}^*(N_v,\widehat{CC}_*(A)).
\eqno(4.33)
$$
Here the direct limit is taken over the iterated cup product 
$$
\cup\alpha:\,{\widehat{\mathbb H}}^*(N_v,\widehat{CC}_*(A))\,\longrightarrow\,{\widehat{\mathbb H}}^*(N_v,\widehat{CC}_*(A))
\eqno(4.34)
$$ 
with the extension class of $\mathrm{(4.10)}$.
\end{theorem}

\begin{proof}

Step 1:\\
Recall the Eilenberg-Zilber theorem in its strong form \cite{Pu} which states that the morphisms (4.21) and (4.22) are in fact natural chain-homotopy equivalences 
if $A$ and $B$ are unital. The chain map $\Delta_{cyc}$ preserves Hodge filtrations and induces on the associated graded Hochschild complexes the chain homotopy equivalence (4.17). As in Step 3 of the proof of Theorem 4.6 one may claim
that this holds for the non-unital algebras $B=k\langle X\rangle,\,X$ a set, as well. In particular, one obtains a natural, continuous, $Z_v$-equivariant chain homotopy equivalence
$$
\widehat{CC}_*(A\langle Z_v\rangle\rtimes U:k\rtimes U)_{\{v\}}\,\overset{\sim}{\longrightarrow}\,\widehat{CC}_*(A)\widehat{\otimes}_k\widehat{CC}_*(k\langle Z_v\rangle\rtimes U:k\rtimes U)_{\{v\}}.
\eqno(4.35)
$$
if $v$ acts trivially on $A$. It follows from Theorem 4.1 that in this case there is a natural isomorphism 
$$
HP_*(A\rtimes G)_{[v]}\,\overset{\simeq}{\longrightarrow}\,\widehat{{\mathbb H}}_*\left(N_v,\widehat{CC}_*(A)\widehat{\otimes}_k \widehat{CC}_*(k\langle Z_v\rangle\rtimes U:k\rtimes U)_{\{v\}} \right).
\eqno(4.36)
$$
Step 2:\\
Consider the composition
$$
\pi_*:\widehat{CC}_*(k\langle Z_v\rangle\rtimes U:k\rtimes U)_{\{v\}}\,\to\,
\widehat{CC}_*(k\langle N_v\rangle\rtimes U:k\rtimes U)_{\{v\}}
$$

$$
\overset{\Delta_{cyc}}{\longrightarrow}\,\widehat{CC}_*(k\langle N_v\rangle)\widehat{\otimes}_k\widehat{CC}(k\rtimes U:k\rtimes U)_{\{v\}}
\,\overset{(2.14)}{\simeq}\,\widehat{CC}_*(k\langle N_v\rangle)
\eqno(4.37)
$$
This $N_v$-equivariant chain map preserves Hodge filtrations because the quasi-isomorphism
$$ 
\widehat{CC}_*(k\langle N_v\rangle:k\langle N_v\rangle)_{\{v\}}\,\overset{\sim}{\longrightarrow}\,k
$$
vanishes on
$\mathrm{Fil}^1_{\mathrm{Hodge}}\widehat{CC}_*(k\langle N_v\rangle:k\langle N_v\rangle)_{\{v\}}$. 
The $N_v$-equivariant continuous linear operator $\pi_*\circ\widetilde{\eta}$
defines then an $N_v$-equivariant chain homotopy between $\pi_*$ and $\cap c_{\alpha}\circ\pi_*$  (see (4.27)), so that altogether we have constructed an $N_v$-equivariant continuous chain map
$$
\Psi:\widehat{CC}_*(k\langle Z_v\rangle\rtimes U:kU)_{\{v\}}\,\longrightarrow\,\mathrm{Cone}(\cap c_{\alpha}-\mathrm{id}:\widehat{CC}_*(k\langle N_v\rangle) \to \widehat{CC}_*(k\langle N_v\rangle)).
\eqno(4.38)
$$
Step 5:\\
Denote by $\widehat{CC}_*(k\langle N_v\rangle)[-2]$ the complex $\widehat{CC}_*(k\langle N_v\rangle)$ with the same grading, but the shifted filtration
$$
\mathrm{Fil}^m_{\mathrm{Hodge}}\widehat{CC}_*(k\langle N_v\rangle)[-2] =\mathrm{Fil}^{m-2}_{\mathrm{Hodge}}\widehat{CC}_*(k\langle N_v\rangle).
$$
This has the effect that $\cap c_{\alpha}:\widehat{CC}_*(k\langle N_v\rangle)\to\widehat{CC}_*(k\langle N_v\rangle)[-2]$ preserves filtrations while 
$\mathrm{id}:\widehat{CC}_*(k\langle N_v\rangle)\to\widehat{CC}_*(k\langle N_v\rangle)[-2]$ shifts them by +2. Consequently
$$
\Psi:\widehat{CC}_*(k\langle Z_v\rangle\rtimes U:kU)_{\{v\}}\,\longrightarrow\,\mathrm{Cone}(\cap c_{\alpha}-\mathrm{id}:\widehat{CC}_*(k\langle N_v\rangle) \to \widehat{CC}_*(k\langle N_v\rangle)[-2])
$$ 
preserves the given filtrations and one may pass to the associated graded situation. Up to natural quasi-isomorphism one finds 
$$
\mathrm{Gr}(\Psi):C_*(k\langle Z_v\rangle\rtimes U:kU)_{\{v\}}\,\longrightarrow\,\mathrm{Cone}(\cap c_{\alpha}:C_*(k\langle N_v\rangle) \to C_*(k\langle N_v\rangle)[-2])
\eqno(4.39)
$$
where $\cap c_{\alpha}$ is the ``classical'' cap-product with the extension cocycle $c_{\alpha}$.\\
Step 6:\\
We recall a classical result of Gysin, which states in our framework that 
$$
\begin{array}{ccccc}
\left(C^{Bar}_*(Z_v,k)\right)_U  & \to & C_*^{Bar}(N_v,k) & \overset{\cap\alpha}{\longrightarrow} & C_*^{Bar}(N_v,k)[-2]  \\
\end{array}
\eqno(4.40)
$$ 
is a distinguished triangle in the chain-homotopy category of (bounded below) complexes of\\ $k\rtimes N_v$-modules. By (4.20) and Lemma 3.2 the same holds for 
$$
\begin{array}{ccccc}
C_*(k\langle Z_v\rangle\rtimes U:kU)_{\{v\}}  & \to & C_*(k\langle N_v\rangle) & \overset{\cap\alpha}{\longrightarrow} & C_*(k\langle N_v\rangle)[-2]. \\
\end{array}
\eqno(4.41)
$$
This is a distinguished triangle of complexes of free $N_v$-modules, so that it remains distinguished after tensoring with $C_*(A)$.
Passing to hyperhomology, one arrives at the vanishing result
$$
{\mathbb H}_*\left(N_v,\mathrm{Cone}(\mathrm{id}_{C_*(A)}\otimes \mathrm{Gr}(\Psi))\right)\,=\,0.
\eqno(4.42)
$$
Iterated use of the five lemma yields then
$$
{\mathbb H}_*(N_v,\mathrm{Cone}(\mathrm{id}_{\widehat{CC}_*(A)}\widehat{\otimes}\Psi)/\mathrm{Fil}^n\mathrm{Cone}(\mathrm{id}_{\widehat{CC}_*(A)}\widehat{\otimes}\Psi))\,=\,0
$$
for $n\geq 0$ and finally 
$$
\widehat{{\mathbb H}}_*(N_v,\mathrm{Cone}(\mathrm{id}_{\widehat{CC}_*(A)}\widehat{\otimes}\Psi))\,=\,0
\eqno(4.43)
$$
by the $lim^1$-exact sequence. \\
Step 7:\\
The vanishing result (4.43), (4.36), and Lemma 3.6 imply a natural isomorphism
$$
HP_*(A\rtimes G)_{[v]}\,\overset{\simeq}{\longrightarrow}\,
\widehat{{\mathbb H}}_*(N_v, \mathrm{Cone}(\cap c_{\alpha}-\mathrm{id}:
\widehat{CC}_*(A)\to \widehat{CC}_*(A))).
\eqno(4.44)
$$
This means that $HP_*(A\rtimes G)_{[v]}$  fits into an exact triangle
$$
\begin{array}{ccc}
\underset{n=0}{\overset{\infty}{\prod}}\,{\mathbb H}_n(N_v,\widehat{CC}_*(A)) &  \overset{\cap\alpha-\mathrm{id}}{\longrightarrow} &  \underset{n=0}{\overset{\infty}{\prod}}\,{\mathbb H}_n(N_v,\widehat{CC}_*(A))\\
 & &  \\
\nwarrow & & \swarrow \delta \\
\end{array}
$$
$$
HP_*(A\rtimes G)_{[v]}
\eqno(4.45)
$$
which is equivalent to the claim
$$
\widehat{{\mathbb H}}_*(N_v, \mathrm{Cone}(\cap c_{\alpha}-\mathrm{id}))\,\overset{\simeq}{\longrightarrow}\,{\mathbb R}\underset{\underset{\cap\alpha}{\longleftarrow}}{\lim}\,\widehat{{\mathbb H}}_*(N_v,\widehat{CC}_*(A)).
\eqno(4.46)
$$
This establishes Theorem 4.8. The proof of Theorem 4.9 is left to the reader.
\end{proof}

\begin{remark}
If one follows the same strategy in the general case one ends up instead at $\mathrm{(4.45)}$ at an exact triangle

$$
\begin{array}{c}
\widehat{{\mathbb H}}_*\left(N_v,{\mathrm Gr}_{0}(\widehat{CC}_*(A\langle Z_v\rangle\rtimes U:k\rtimes U)_{\{v\}})\right) \,\overset{\delta}{\longrightarrow} \,  
\widehat{{\mathbb H}}_*\left(N_v,{\mathrm Gr}_{1}(\widehat{CC}_*(A\langle Z_v\rangle\rtimes U:k\rtimes U)_{\{v\}})\right)\\
    \\
\nwarrow \,\,\,\,\,\,\,\,\,\,\,\,\,\,\,\,\,\,\,\,\,\,\,\,\,\,\,\,\,\,\,\,\,\,\,\,\,\,\,\,\,\,\,\,\,\,\,\,\,\,\,\,\,\,\,\,\,\,\,\,\,\,\,\,\,\,\,\,\,\,\,\,\,\,\,\,\,\,\,\,\,\,\,\, \swarrow\delta \\
 HP_*(A\rtimes G)_{[v]}\\
\end{array}
\eqno(4.47)
$$

\end{remark}


\begin{thebibliography}{Wittgenstein} 
\newcommand{\hotz}[1]{\bibitem[{#1}]{#1}} 

\hotz{Bu} {\sc D.~Burghelea,} The cyclic homology of group rings,\\ Comment. Math. Helv. 60 (1985), 354--365

\hotz{Co} {\sc A.~Connes,} Noncommutative differential geometry,\\ 
 Publ. Math. IHES 62, (1985), 41--144 

\hotz{CQ1} {\sc J.~Cuntz, D.~Quillen,} Algebra extensions and nonsingularity,\\ J.AMS 8 (1995), 251--289

\hotz{CQ2}  {\sc J.~Cuntz, D.~Quillen,} Cyclic homology and nonsingularity,\\ J.AMS 8 (1995), 373--442

\hotz{KS} {\sc M.~Kashiwara, P.~Shapira}, Sheaves on manifolds,\\ Springer Grundlehren 292 (1990)

\hotz{Lo} {\sc J.L.~Loday}, Cyclic homology,\\ Springer Grundlehren 301 (1992)

\hotz{Ni}  {\sc V.~Nistor,} Group cohomology and the cyclic
cohomology of crossed products, Invent. Math. 99,  (1990), 411--423

\hotz{Pu} {\sc M.~Puschnigg}, Explicit product structures in cyclic homology theories, J.K-Theory 15 (1998), 323--345


\end{thebibliography}
\end{document}